\documentclass[a4paper,11pt]{article}
\usepackage{amsmath}
\usepackage{amsfonts}

\input{tcilatex}
\begin{document}

\begin{center}
$f-$\textbf{Statistical Convergence of order }$\beta $

\textbf{for Sequences of Fuzzy Numbers}

\textit{\bigskip}

{\footnotesize H\i fs\i \ ALTINOK}$^{1}${\footnotesize \qquad Mithat KASAP}$%
^{2}$

$^{1}${\footnotesize Department of Mathematics, F\i rat University, 23119,
Elaz\i \u{g}, TURKEY.}

$^{2}${\footnotesize Department of Accounting, \c{S}\i rnak University,
23119, \c{S}\i rnak, TURKEY.}

{\footnotesize E-mail: hifsialtinok@yahoo.com\quad fdd\_mithat@hotmail.com}

\bigskip
\end{center}

\begin{quote}
\textbf{Abstract.} In this paper, we extend the notions of statistically
convergence of order $\beta $ and strong Ces\`{a}ro summability of order $%
\beta ,$ and introduce the notions $f-$statistically convergence of order $%
\beta $ and strong Ces\`{a}ro summability of order $\beta $ for $\beta \in
\left( 0,1\right] $ with respect to an unbounded modulus function $f$ for
sequences of fuzzy numbers and give some inclusion theorems.
\end{quote}

\bigskip

\noindent \textbf{Key words and phrases:} Fuzzy number, sequence of fuzzy
numbers, statistical convergence, Ces\`{a}ro summability, modulus function.

\noindent \textbf{Mathematics Subject Classification:} 40A05; 40A25; 40A30;
40C05; 03E72.

\bigskip

\noindent \textbf{1. Introduction}

$\medskip $

In order to generalize the concept of convergence of real sequences, the
notion of statistical convergence was introduced by Fast \cite{fast} and
Schoenberg \cite{schonberg}, independently. Over the years and under
different names statistical convergence has been discussed in the theory of
Fourier analysis, Ergodic theory and Number theory. Later on it was further
investigated from the sequence space point of view and linked with
summability theory by Fridy \cite{fridy}, \v{S}al\'{a}t \cite{Salat}, Connor 
\cite{connor}, topological groups by \c{C}akall\i \  \cite{cakalli}, function
spaces by Caserta \textit{et al.} \cite{caserta}.

Matloka \cite{matloka} defined the notion of fuzzy sequence and introduced
bounded and convergent sequences of fuzzy real numbers and studied their
some properties. After then, Nuray and Sava\c{s} \cite{Nuray and savas}
defined the notion of statistical convergence for sequences of fuzzy
numbers. Since then, there has been increasing interest in the study of
statistical convergence of fuzzy sequences (see \cite{tripathy1},\cite%
{canak1},\cite{mursaleen2},\cite{savas},\cite{hazarika1},\cite{aytar}).

\bigskip

\c{C}olak \cite{colak} generalized the statistical convergence by ordering
the interval $\left[ 0,1\right] $ and defined the statistical convergence of
order $\alpha $ and strong $p-$Ces\`{a}ro summability of order $\alpha $,
where $0<\alpha \leq 1$ and $p$ is a positive real number. Altinok \textit{%
et al.} \cite{altinok1} introduced the concepts of statistical convergence
of order $\beta $ and strong $p-$Ces\`{a}ro summability of order $\beta $\
for sequences of fuzzy numbers. Aizpuru \textit{et al. }\cite{aizpuru}%
\textit{\ }defined the $f-$density of the subset $A$ of $%
\mathbb{N}
$ by using an unbounded modulus function. After then, Bhardwaj \cite%
{bhardwaj} introduced $f-$statistical convergence of order $\alpha $ and
strong Cesaro summability of order $\alpha $ with respect to a modulus
function $f$ for real sequences. The purpose of this paper is to generalize
the study of Bhardwaj \cite{bhardwaj} and \c{C}olak \cite{colak} applying to
sequences of fuzzy numbers so as to fill up the existing gaps in the
summability theory of fuzzy numbers. For a detailed account of many more
interesting investigations concerning statistical convergence of order $%
\alpha $ and $\beta $, one may refer to (\cite{colak2},\cite{et1},\cite{et2},%
\cite{et3},\cite{altin1},\cite{altin2},\cite{altinok2})

\bigskip

This paper organizes as follows: In section 2, we give the basic notions
which will be used throughout the paper. In section 3, we define the spaces $%
S^{\beta }\left( F,f\right) $ and $w^{\beta }\left( F,f\right) ,$ the set of
all $f-$statistically convergent sequences of order $\beta $ and strong Ces%
\`{a}ro summability of order $\beta $ with respect to an unbounded modulus
function $f$ for fuzzy sequences, respectively and establish inclusion
relations among the spaces\ $S^{\beta }\left( F,f\right) $ for different
values of $\beta $. Moreover, we obtain inclusion relations between the
spaces $w^{\beta }\left( F,f\right) $, $w^{\beta ,0}\left( F,f\right) $ and $%
w^{\beta ,\infty }\left( F,f\right) $\ and give some conditions related to
modulus function $f$ for inclusion relation $w^{\beta }\left( F,f\right)
\subset S^{\gamma }\left( F,f\right) .$

\bigskip

\noindent \textbf{2. Definitions and Preliminaries}

\bigskip

In this section, we recall some basic definitions and notations that we are
going to use in this paper.

\bigskip

A fuzzy set $u$ on $\mathbb{R}$ is called a fuzzy number if it has the
following properties:

$i)$ $u$ is normal, that is, there exists an $x_{0}\in \mathbb{R}$ such that 
$u(x_{0})=1;$

$ii)$ $u$ is fuzzy convex, that is, for $x,y\in \mathbb{R}$ and $0\leq
\lambda \leq1,u(\lambda x+(1-\lambda)y)\geq \min[u(x),u(y)];$

$iii)$ $u$ is upper semicontinuous;

$iv)$ $\limfunc{supp}u=cl\{x\in \mathbb{R}:u(x)>0\},$ or denoted by $%
[u]^{0}, $ is compact.

$\alpha-$level set $[u]^{\alpha}$ of a fuzzy number $u$ is defined by 
\begin{equation*}
\lbrack u]^{\alpha}=\left \{ 
\begin{array}{cc}
\{x\in \mathbb{R}:u(x)\geq \alpha \}, & \text{if }\alpha \in \left( 0,1%
\right] \\ 
\limfunc{supp}u, & \text{if }\alpha=0.%
\end{array}
\right.
\end{equation*}
It is clear that $u$ is a fuzzy number if and only if $[u]^{\alpha}$ is a
closed interval for each $\alpha \in \left[ 0,1\right] $ and $[u]^{1}\neq
\emptyset.$\ We denote space of all fuzzy numbers by $L(\mathbb{R}).$

In order to calculate the distance between two fuzzy numbers $u$ and $v,$ we
use the metric%
\begin{equation*}
d\left( u,v\right) =\underset{0\leq \alpha \leq1}{\sup}d_{H}\left( \left[ u%
\right] ^{\alpha},\left[ v\right] ^{\alpha}\right)
\end{equation*}
where $d_{H}$ is the Hausdorff metric defined by%
\begin{equation*}
d_{H}\left( \left[ u\right] ^{\alpha},\left[ v\right] ^{\alpha}\right) =\max
\left \{ \left \vert \underline{u}^{\alpha}-\underline{v}^{\alpha }\right
\vert ,\left \vert \overline{u}^{\alpha}-\overline{v}^{\alpha }\right \vert
\right \} .
\end{equation*}
It is known that $d$ is a metric on $L(\mathbb{R}),$ and $\left( L(\mathbb{R}%
),d\right) $ is a complete metric space.

A sequence $X=(X_{k})$ of fuzzy numbers is a function $X:\mathbb{%
N\rightarrow }L(\mathbb{R}).$ Let $X=\left( X_{k}\right) $ be a sequence of
fuzzy numbers. The sequence $X=\left( X_{k}\right) $ of fuzzy numbers is
said to be bounded if the set $\left \{ X_{k}:k\in \mathbb{N}\right \} $ of
fuzzy numbers is bounded and convergent to the fuzzy number $X_{0}$, written
as $\lim \limits_{k}X_{k}$ $=X_{0}$, if for every $\varepsilon >0$ there
exists a positive integer $k_{0}$ such that $d\left( X_{k},X_{0}\right)
<\varepsilon $ for $k>k_{0}.$ Let $s\left( F\right) ,$ $\ell _{\infty
}\left( F\right) $ and $c\left( F\right) $ denote the set of all sequences,
all bounded sequences and all convergent sequences of fuzzy numbers,
respectively \cite{matloka}.

\bigskip

The concept of modulus function was formally introduced by Nakano \cite%
{nakano}. A mapping $f:[0,\infty)$ $\rightarrow$ $[0,\infty)$ is said to be
a modulus if

$i)$ $f\left( x\right) =0$ iff $x=0$,

$ii)$ $f\left( x+y\right) \leq f\left( x\right) +f\left( y\right) $ for $%
x,y\geq0,$

$iii)$ $f$ is increasing,

$iv)$ $f$ is right-continuous at $0.$

The continuity of $f$ everywhere on $[0,\infty )$ follows from above
definition. A modulus function can be bounded or unbounded. For example $%
f\left( x\right) =x^{p},\left( 0<p\leq 1\right) $ is bounded and $f\left(
x\right) =\frac{x}{1+x}$ is bounded. For an extensive view on this subject
we refer (\cite{altin1},\cite{Sarma},\cite{talo2},\cite{et1},\cite{et2}).

\bigskip

We recall notions of statistically convergence of order $\beta $ and $%
f_{\beta }-$density.

\bigskip

\noindent \textbf{Definition 2.1.} \cite{altinok1} Let $\beta \in \left( 0,1%
\right] $ and $X=\left( X_{k}\right) $ be a sequence of fuzzy numbers. Then
the sequence $X=\left( X_{k}\right) $ of fuzzy numbers is said to be
statistically convergent of order $\beta ,$ to fuzzy number $X_{0}$ if for
every $\varepsilon >0,$%
\begin{equation*}
\lim_{n\rightarrow \infty }\frac{1}{n^{\beta }}\left \vert \left \{ k\leq
n:d\left( X_{k},X_{0}\right) \geq \varepsilon \right \} \right \vert =0,
\end{equation*}%
where the vertical bars indicate the number of elements in the enclosed set.
In this case we write $S^{\beta }\left( F\right) -\lim X_{k}=X_{0}.$ We
denote the set of all statistically convergent sequences of order $\beta $
by $S^{\beta }\left( F\right) .$

\bigskip

We take $\beta $ instead of $\alpha $ for suitability with our spaces.

\bigskip

\noindent \textbf{Definition 2.2. }\cite{bhardwaj}\textbf{\ }Let $f$ be an
unbounded modulus function and $\beta \in \left( 0,1\right] $ be any real
number. The $f_{\beta }-$density of the subset $A$ of $%
\mathbb{N}
$ defined by

\begin{equation*}
d_{\beta }^{f}\left( A\right) =\lim_{n\rightarrow \infty }\frac{1}{f\left(
n^{\beta }\right) }f\left( \left \vert \left \{ k\leq n:k\in A\right \}
\right \vert \right) .
\end{equation*}

\noindent \textbf{Remark 2.3. }\cite{bhardwaj}\textbf{\ }The $f_{\beta }-$%
density reduces to the natural density for $\beta =1$ and $f\left( x\right)
=x$, and $f_{\beta }-$density becomes the $f-$density in case $\beta =1$. $%
f_{\beta }-$density is $\beta -$density when $f\left( x\right) =x.$

Generally, the equality $d_{\beta }^{f}\left( A\right) +d_{\beta }^{f}\left( 
\mathbb{N}
-A\right) =1$ does not hold for any unbounded modulus $f$. For example, if
we take $f\left( x\right) =x^{p}$ for $0<p\leq 1$, $\beta \in \left(
0,1\right) $ and $A=\left \{ 2n:n\in 
\mathbb{N}
\right \} $, then $d_{\beta }^{f}\left( A\right) =\infty =d_{\beta
}^{f}\left( 
\mathbb{N}
-A\right) $. Moreover, finite sets have zero $f_{\beta }-$density for any
unbounded modulus function $f.$

\bigskip

\noindent \textbf{Remark 2.4. }\cite{bhardwaj}\textbf{\ }If $d_{\beta
}^{f}\left( A\right) =0,$ then $d_{\beta }\left( A\right) =0$, and so $%
d\left( A\right) =0,$ where $f$ is any unbounded modulus function and $\beta
\in \left( 0,1\right] $. Conversely, a set having zero natural density may
have non-zero $f_{\beta }-$density for same $f$ and $\beta $. For this, we
give following example.

\bigskip

\noindent \textbf{Example 2.5.} We consider the modulus function $f\left(
x\right) =\ln \left( x+1\right) $ and the set $A=\left \{
1,8,27,64,...\right \} $. Then, it is easy to see that $d\left( A\right) =0$
and $d_{\beta }\left( A\right) =0$ for $\beta \in \left( \frac{1}{3},1\right]
,$ but $d_{\beta }^{f}\left( A\right) \geq d^{f}\left( A\right) =\frac{1}{3}%
, $ so $d_{\beta }^{f}\left( A\right) \neq 0$.

\bigskip

\noindent \textbf{Lemma 2.6.} \cite{bhardwaj}\textbf{\ }Let\textbf{\ }$f$ be
any unbounded modulus, $0<\beta \leq \gamma \leq 1$ and $A\subset 
\mathbb{N}
$. Then $d_{\gamma }^{f}\left( A\right) \leq d_{\beta }^{f}\left( A\right) $.

Thus, if the set $A$ has zero $f_{\beta }-$density for any unbounded modulus 
$f$ and $0<\beta \leq \gamma \leq 1$, then it has zero $f_{\gamma }-$%
density. In particular, for some $\beta \in \left( 0,1\right] ,$ a set
having zero $f_{\beta }-$density has zero $f-$density. But, the converse is
not true.

\bigskip

\noindent \textbf{3. Main Results}

\bigskip

We now introduce the concept of $f-$statistical convergence of order $\beta $
for sequences of fuzzy numbers as follows.

\bigskip

\noindent \textbf{Definition 3.1.} Let $X=\left( X_{k}\right) $ be a fuzzy
sequence, $f$ be an unbounded modulus and $\beta \in \left( 0,1\right] $. A
sequence $X=\left( X_{k}\right) $ of fuzzy numbers is said to be $f-$%
statistically convergent of order $\beta $ to fuzzy number $X_{0}$ or $%
S^{\beta }\left( F,f\right) -$convergent to $X_{0}$ if for each $\varepsilon
>0$,%
\begin{equation*}
\lim_{n\rightarrow \infty }\frac{1}{f\left( n^{\beta }\right) }f\left(
\left \vert \left \{ k\leq n:d\left( X_{k},X_{0}\right) \geq \varepsilon
\right \} \right \vert \right) =0.
\end{equation*}%
In this case, we write $S^{\beta }\left( F,f\right) -\lim X_{k}=X_{0}$. By $%
S^{\beta }\left( F,f\right) ,$ we shall denote the set of all sequences of
fuzzy numbers which are $f-$statistically convergent of order $\beta $ for
convenience with our previous studies and by $S^{\beta ,0}\left( F,f\right) ,
$ the set of all $f-$statistically null sequences of order $\beta $. For any
unbounded modulus function $f$ and $\beta \in \left( 0,1\right] ,$ the
inclusion relation $S^{\beta ,0}\left( F,f\right) \subset S^{\beta }\left(
F,f\right) $ is clear. Furthermore, we point out $f-$statistical convergence
of order $\beta $ reduces statistical convergence of order $\beta $ defined
in \cite{altinok1} for $f\left( x\right) =x$ and $f-$statistical convergence
of order $\beta $ reduces to statistical convergence defined in \cite{Nuray
and savas} for $f\left( x\right) =x$ and $\beta =1$.

\bigskip

Now, we generalize the space of strongly Ces\`{a}ro summable sequences of
order $\beta $ using a modulus function.

\bigskip

\noindent \textbf{Definition 3.2.} Let $X=\left( X_{k}\right) $ be a
sequence of fuzzy numbers, $f$ be an unbounded modulus function and $\beta
\in \left( 0,1\right] $. We define following spaces: 
\begin{eqnarray*}
w^{\beta ,0}\left( F,f\right) &=&\left \{ X\in s\left( F\right)
:\lim_{n\rightarrow \infty }\frac{1}{n^{\beta }}\dsum
\limits_{k=1}^{n}f\left( d\left( X_{k},\bar{0}\right) \right) =0\right \} ,
\\
w^{\beta }\left( F,f\right) &=&\left \{ X\in s\left( F\right)
:\lim_{n\rightarrow \infty }\frac{1}{n^{\beta }}\dsum
\limits_{k=1}^{n}f\left( d\left( X_{k},X_{0}\right) \right) =0\text{ for
some }X_{0}\in L\left( \mathbb{R}\right) \right \} , \\
w^{\beta ,\infty }\left( F,f\right) &=&\left \{ X\in s\left( F\right)
:\sup_{n}\frac{1}{n^{\beta }}\dsum \limits_{k=1}^{n}f\left( d\left( X_{k},%
\bar{0}\right) \right) <\infty \right \} ,
\end{eqnarray*}%
In case $f\left( x\right) =x$, the space $w^{\beta }\left( F,f\right) $ is
the same as the space $w^{\beta }\left( F,p\right) $ of Altinok \textit{et
al.} \cite{altinok1} for $p=1$. We shall denote the space $w^{\beta ,\infty
}\left( F,f\right) $ by $w^{\beta ,\infty }\left( F\right) $ for $f\left(
x\right) =x$. It is easy to see that $w^{\beta }\left( F,f\right) $, $%
w^{\beta ,0}\left( F,f\right) $ and $w^{\beta ,\infty }\left( F,f\right) $
are linear spaces over the complex field $%
\mathbb{C}
$.

\bigskip

\noindent \textbf{Remark 3.3.} In Definition 3.1, although the $f-$%
statistical convergence of order $\beta $ is well defined for $\beta \in
\left( 0,1\right] ,$ but not well defined for $\beta >1.$ (See Example 3.4).

\bigskip

\noindent \textbf{Example 3.4.} Let $f$\ be an unbounded modulus function
such that $\lim \limits_{t\rightarrow \infty }\frac{f\left( t\right) }{t}>0$%
\ and $X=\left( X_{k}\right) $\ be a sequence of fuzzy numbers as follows:%
\begin{equation*}
X_{k}\left( x\right) =\left \{ 
\begin{array}{cc}
\begin{array}{c}
x+3, \\ 
-x-1, \\ 
0%
\end{array}
& \left. 
\begin{array}{c}
-3\leq x\leq -2 \\ 
-2\leq x\leq -1 \\ 
\text{otherwise}%
\end{array}%
\right \} :X_{0}\text{, if }k\text{ is odd} \\ 
\begin{array}{c}
x-1, \\ 
-x+3, \\ 
0%
\end{array}
& \left. 
\begin{array}{c}
1\leq x\leq 2 \\ 
2\leq x\leq 3 \\ 
\text{otherwise}%
\end{array}%
\right \} :X_{0}^{\prime }\text{, if }k\text{ is even}%
\end{array}%
\right.
\end{equation*}%
If we calculate the $\alpha -$level set of sequence $\left( X_{k}\right) ,$
then we find the set%
\begin{equation*}
\left[ X_{k}\right] ^{\alpha }=\left \{ 
\begin{array}{cc}
\left[ \alpha -3,-1-\alpha \right] : & \left[ X_{0}\right] ^{\alpha }\text{,
if }k\text{ is odd} \\ 
\left[ \alpha +1,3-\alpha \right] : & \left[ X_{0}^{\prime }\right] ^{\alpha
}\text{, if }k\text{ is even}%
\end{array}%
\right.
\end{equation*}%
Then we can write 
\begin{equation*}
\frac{1}{f\left( n^{\beta }\right) }f\left( \left \vert \left \{ k\leq
n:d\left( X_{k},X_{0})\geq \varepsilon \right) \right \} \right \vert
\right) \leq \frac{f\left( \frac{n}{2}\right) }{f\left( n^{\beta }\right) }
\end{equation*}%
and 
\begin{equation*}
\frac{1}{f\left( n^{\beta }\right) }f\left( \left \vert \left \{ k\leq
n:d\left( X_{k},X_{0}^{\prime })\geq \varepsilon \right) \right \} \right
\vert \right) \leq \frac{f\left( \frac{n}{2}\right) }{f\left( n^{\beta
}\right) }
\end{equation*}%
So, we get%
\begin{equation*}
\lim_{n\rightarrow \infty }\frac{1}{f\left( n^{\beta }\right) }f\left( \left
\vert \left \{ k\leq n:d\left( X_{k},X_{0})\geq \varepsilon \right) \right
\} \right \vert \right) =0
\end{equation*}%
and%
\begin{equation*}
\lim_{n\rightarrow \infty }\frac{1}{f\left( n^{\beta }\right) }f\left( \left
\vert \left \{ k\leq n:d\left( X_{k},X_{0}^{\prime })\geq \varepsilon
\right) \right \} \right \vert \right) =0
\end{equation*}%
for $\beta >1$ and for each $\varepsilon >0$ using the property $\lim
\limits_{t\rightarrow \infty }\frac{f\left( t\right) }{t}>0.$ Therefore,
sequence $X=\left( X_{k}\right) $ of fuzzy numbers is $f-$statistically
convergent of order $\beta $ to both $X_{0}$ and $X_{0}^{\prime }.$ That is,
the $f-$statistical limit of order $\beta $ may not be unique for $\beta >1$
(See Fig. 3.1). 
\begin{equation*}
\FRAME{itbpF}{4.5343in}{2.1651in}{0in}{}{}{Figure}{\special{language
"Scientific Word";type "GRAPHIC";maintain-aspect-ratio TRUE;display
"USEDEF";valid_file "T";width 4.5343in;height 2.1651in;depth
0in;original-width 3.7182in;original-height 1.7598in;cropleft "0";croptop
"1";cropright "1";cropbottom "0";tempfilename
'OBFPOP02.wmf';tempfile-properties "XPR";}}
\end{equation*}

Every convergent sequence of fuzzy numbers is $f-$statistically convergent
of order $\beta $ for any unbounded modulus function $f$ and $0<\beta \leq
1. $ But, the converse is not true, for this we can give the following
example.

\bigskip

\noindent \textbf{Example 3.5.} Define the sequence $X=\left( X_{k}\right) $
of fuzzy numbers as follows:%
\begin{equation*}
X_{k}\left( x\right) =\left \{ 
\begin{array}{cc}
\begin{array}{c}
x+3, \\ 
-x-1, \\ 
0%
\end{array}
& \left. 
\begin{array}{c}
-3\leq x\leq -2 \\ 
-2\leq x\leq -1 \\ 
\text{otherwise}%
\end{array}%
\right \} \text{ if }k=n^{3} \\ 
\begin{array}{c}
x-1, \\ 
-x+3, \\ 
0%
\end{array}
& \left. 
\begin{array}{c}
1\leq x\leq 2 \\ 
2\leq x\leq 3 \\ 
\text{otherwise}%
\end{array}%
\right \} \text{ if }k\neq n^{3}%
\end{array}%
\right. 
\end{equation*}%
Take modulus function $f\left( x\right) =x^{p}$ for $0<p\leq 1.$ We can find
the $\alpha -$level set of fuzzy sequence $\left( X_{k}\right) $ as follows:%
\begin{equation*}
\left[ X_{k}\right] ^{\alpha }=\left \{ 
\begin{array}{cc}
\left[ \alpha -3,-1-\alpha \right] , & \text{if }k=n^{3}\text{ } \\ 
\left[ \alpha +1,3-\alpha \right] , & \text{if }k\neq n^{3}\text{ }%
\end{array}%
\right. 
\end{equation*}%
Then, the fuzzy sequence $\left( X_{k}\right) $ is $f-$statistically
convergent of order $\beta $ for $\beta \in \left( \frac{1}{3},1\right] ,$
but not convergent.

\bigskip

\noindent \textbf{Theorem 3.6.} Let $X=\left( X_{k}\right) $, $Y=\left(
Y_{k}\right) $ be any two fuzzy sequences and $L_{1},$ $L_{2}$ be fuzzy
numbers. Also, let $f$ be an unbounded modulus function and $\beta \in
\left( 0,1\right] $. Then

$\left( i\right) $ If $S^{\beta }\left( F,f\right) -\lim X_{k}=L_{1}$ and $%
c\in 
\mathbb{C}
$, then $S^{\beta }\left( F,f\right) -\lim cX_{k}=cL_{1}$.

$\left( ii\right) $ If $S^{\beta }\left( F,f\right) -\lim X_{k}=L_{1}$ and $%
S^{\beta }\left( F,f\right) -\lim Y_{k}=L_{2}$, then $S^{\beta }\left(
F,f\right) -\lim \left( X_{k}+Y_{k}\right) =L_{1}+L_{2}$.

\bigskip

\noindent \textbf{Theorem 3.7.} Let $f$ be unbounded modulus function and $%
\beta ,$ $\gamma $ be real numbers such that $0<\beta \leq \gamma \leq 1$.
Then $S^{\beta }\left( F,f\right) \subset S^{\gamma }\left( F,f\right) $ and
the inclusion is strict.

\textbf{Proof.} It can be easily shown the inclusion by using the fact that $%
f$ is increasing for $0<\beta \leq \gamma \leq 1$. Now, we show that the
inclusion is strict. For this, consider fuzzy sequence $X=\left(
X_{k}\right) $ defined by

\begin{equation*}
X_{k}\left( x\right) =\left \{ 
\begin{array}{cc}
\begin{array}{c}
x+3, \\ 
-x-1, \\ 
0%
\end{array}
& \left. 
\begin{array}{c}
-3\leq x\leq -2 \\ 
-2\leq x\leq -1 \\ 
\text{otherwise}%
\end{array}%
\right \} \text{ if }k=n^{2} \\ 
\begin{array}{c}
x-1, \\ 
-x+3, \\ 
0%
\end{array}
& \left. 
\begin{array}{c}
1\leq x\leq 2 \\ 
2\leq x\leq 3 \\ 
\text{otherwise}%
\end{array}%
\right \} \text{ if }k\neq n^{2}%
\end{array}%
\right. 
\end{equation*}%
and take modulus function $f\left( x\right) =x^{p}$, $0<p\leq 1$. We can
find the $\alpha -$level set of sequence $\left( X_{k}\right) $ as follows: 
\begin{equation*}
\left[ X_{k}\right] ^{\alpha }=\left \{ 
\begin{array}{cc}
\left[ \alpha -3,-1-\alpha \right] , & \text{if }k=n^{2} \\ 
\left[ \alpha +1,3-\alpha \right] , & \text{if }k\neq n^{2}%
\end{array}%
\right. 
\end{equation*}%
Then, the fuzzy sequence $\left( X_{k}\right) $ is $f-$statistically
convergent of order $\gamma $ for $\gamma \in \left( \frac{1}{2},1\right] ,$
but not $f-$statistically convergent of order $\beta $ for $\beta \in \left(
0,\frac{1}{2}\right] .$

\bigskip

\noindent \textbf{Corollary 3.8.} Let $X=\left( X_{k}\right) $ be a fuzzy
sequence, $f$ be an unbounded modulus function and $\beta \in \left( 0,1%
\right] $. Then $S^{\beta }\left( F,f\right) \subset S\left( F,f\right) $
and the inclusion is strict, also the limits of sequence $X=\left(
X_{k}\right) $ of fuzzy numbers are same.

\bigskip

We have the following theorem from Remark 2.4.

\bigskip

\noindent \textbf{Theorem 3.9.} Let $f$ be an unbounded modulus function and 
$\beta \in \left( 0,1\right] $. Then

$\left( i\right) $ $S^{\beta }\left( F,f\right) \subset S^{\beta }\left(
F\right) $ and the inclusion is strict.

$\left( ii\right) $ $S^{\beta }\left( F,f\right) \subset S\left( F\right) $
and the inclusion is strict.

\textbf{Proof.} To show that the strictness of inclusion, consider the fuzzy
sequence $X=\left( X_{k}\right) $ defined as follows%
\begin{equation*}
X_{k}\left( x\right) =\left \{ 
\begin{array}{cc}
\left. 
\begin{array}{cc}
x-k+1, & k-1\leq x\leq k \\ 
-x+k+1, & k\leq x\leq k+1 \\ 
0, & d.d%
\end{array}%
\right \}  & \text{if }k=n^{3} \\ 
\bar{0} & \text{if }k\neq n^{3}%
\end{array}%
\right. 
\end{equation*}%
The $\alpha -$level set of sequence $\left( X_{k}\right) $ is%
\begin{equation*}
\left[ X_{k}\right] ^{\alpha }=\left \{ 
\begin{array}{cc}
\left[ \alpha +k-1,k+1-\alpha \right] , & \text{if }k=n^{3} \\ 
\bar{0}, & \text{otherwise}%
\end{array}%
\right. 
\end{equation*}%
Take modulus function $f\left( x\right) =\ln \left( x+1\right) .$ Then we
see that the sequence $X=\left( X_{k}\right) $ is statistically convergent
of order $\beta $ for $\beta \in \left( \frac{1}{3},1\right] $ and so it is
statistically convergent (See Fig. 3.2). However, $X=\left( X_{k}\right) $
is not $f-$statistically convergent of order $\beta $ since 
\begin{eqnarray*}
d_{\beta }^{f}\left( \left \{ k\in 
\mathbb{N}
:d\left( X_{k},\bar{0}\right) \geq \varepsilon \right \} \right)  &\geq
&d^{f}\left( \left \{ k\in 
\mathbb{N}
:d\left( X_{k},\bar{0}\right) \geq \varepsilon \right \} \right)  \\
&=&\frac{1}{3}\neq 0.
\end{eqnarray*}%
\begin{equation*}
\FRAME{itbpF}{5.3299in}{1.9567in}{0in}{}{}{Figure}{\special{language
"Scientific Word";type "GRAPHIC";maintain-aspect-ratio TRUE;display
"USEDEF";valid_file "T";width 5.3299in;height 1.9567in;depth
0in;original-width 5.4035in;original-height 1.9655in;cropleft "0";croptop
"1";cropright "1";cropbottom "0";tempfilename
'OBFPI101.wmf';tempfile-properties "XPR";}}
\end{equation*}

Now, we give some results related to spaces $w^{\beta }\left( F,f\right) $, $%
w^{\beta ,0}\left( F,f\right) $ and $w^{\beta ,\infty }\left( F,f\right) $
introduced in Definition 3.2.

\bigskip

\noindent \textbf{Remark 3.10.} We didn't allowed $\beta $ to exceed $1$ in
the $w^{\beta }\left( F,p\right) $ of Altinok \textit{et al}. \cite{altinok1}%
, but we consider $\beta $ as any positive real number and it can exceed $1$
in the spaces $w^{\beta ,0}\left( F,f\right) $ and $w^{\beta }\left(
F,f\right) $.

\bigskip

\noindent \textbf{Theorem 3.11. }Let $f$ be any modulus function. Then

$\left( i\right) $ $w^{\beta ,0}\left( F,f\right) \subset w^{\beta ,\infty
}\left( F,f\right) $ for $\beta >0,$

$\left( ii\right) $ $w^{\beta }\left( F,f\right) \subset w^{\beta ,\infty
}\left( F,f\right) $ for $\beta \geq 1.$

\textbf{Proof.} The proof of $(i)$ is trivial, so we give the proof of $(ii)$%
. For this, take $\beta \geq 1$ and any fuzzy sequence $X=\left(
X_{k}\right) $ in the space $w^{\beta }\left( F,f\right) $. Then, we obtain

\begin{equation*}
\frac{1}{n^{\beta }}\dsum \limits_{k=1}^{n}f\left( d\left( X_{k},\bar{0}%
\right) \right) \leq \frac{1}{n^{\beta }}\dsum \limits_{k=1}^{n}f\left(
d\left( X_{k},X_{0}\right) \right) +f\left( d\left( X_{0},\bar{0}\right)
\right) \frac{1}{n^{\beta }}\dsum \limits_{k=1}^{n}1,
\end{equation*}%
from definition of modulus function, and so we have $X\in w^{\beta ,\infty
}\left( F,f\right) $.

\bigskip

\noindent \textbf{Theorem 3.12.} Let $f$ be any modulus function and $\beta
\geq 1$. Then, we have the inclusion relations $w^{\beta }\left( F\right)
\subset w^{\beta }\left( F,f\right) $, $w^{\beta ,0}\left( F\right) \subset
w^{\beta ,0}\left( F,f\right) $ and $w^{\beta ,\infty }\left( F\right)
\subset w^{\beta ,\infty }\left( F,f\right) $.

\textbf{Proof.} We shall prove the inclusion $w^{\beta ,\infty }\left(
F\right) \subset w^{\beta ,\infty }\left( F,f\right) $ since the proofs of
first two inclusion relations are easy. For this, take $\beta \geq 1$ and
any fuzzy sequence $X=\left( X_{k}\right) $ in the space $w^{\beta ,\infty
}\left( F\right) $ so that%
\begin{equation*}
\sup_{n}\frac{1}{n^{\beta }}\dsum \limits_{k=1}^{n}d\left( X_{k},\bar{0}%
\right) <\infty 
\end{equation*}%
Given $\varepsilon >0$ and $\delta \in \left( 0,1\right) $ such that $%
f\left( t\right) <\varepsilon $ for $t\in \left( 0,\delta \right] $.
Consider 
\begin{equation*}
\dfrac{1}{n^{\beta }}\dsum \limits_{k=1}^{n}f\left( d\left( X_{k},\bar{0}%
\right) \right) =\dsum \limits_{1}+\dsum \limits_{2},
\end{equation*}%
where $\dsum \limits_{1}$ is over $d\left( X_{k},\bar{0}\right) \leq \delta $
and $\dsum \limits_{2}$ is over $d\left( X_{k},\bar{0}\right) >\delta $. Then 
$\dsum \limits_{1}\leq \varepsilon \dfrac{1}{n^{\beta -1}}$ and we can write%
\begin{equation*}
d\left( X_{k},\bar{0}\right) <\frac{d\left( X_{k},\bar{0}\right) }{\delta }%
<1+\left[ \left \vert \frac{d\left( X_{k},\bar{0}\right) }{\delta }%
\right \vert \right] ,
\end{equation*}%
for $d\left( X_{k},\bar{0}\right) >\delta ,$ where $\left[ \left \vert
d\left( X_{k},\bar{0}\right) /\delta \right \vert \right] $ denotes the
integer part of $d\left( X_{k},\bar{0}\right) /\delta $. Therefore, for $%
d\left( X_{k},\bar{0}\right) >\delta ,$ we obtain%
\begin{equation*}
f\left( d\left( X_{k},\bar{0}\right) \right) \leq \left( 1+\left[ \left \vert 
\frac{d\left( X_{k},\bar{0}\right) }{\delta }\right \vert \right] \right)
f\left( 1\right) \leq 2f\left( 1\right) \frac{d\left( X_{k},\bar{0}\right) }{%
\delta }
\end{equation*}%
from definition of modulus function. Hence we get inequality 
\begin{equation*}
\dsum \limits_{2}\leq 2f\left( 1\right) \delta ^{-1}\dfrac{1}{n^{\beta }}%
\dsum \limits_{k=1}^{n}d\left( X_{k},\bar{0}\right) 
\end{equation*}%
which together with $\dsum \limits_{1}\leq \varepsilon \dfrac{1}{n^{\beta -1}}
$ yields%
\begin{equation*}
\frac{1}{n^{\beta }}\dsum \limits_{k=1}^{n}f\left( d\left( X_{k},\bar{0}%
\right) \right) \leq \varepsilon \frac{1}{n^{\beta -1}}+2f\left( 1\right)
\delta ^{-1}\frac{1}{n^{\beta }}\dsum \limits_{k=1}^{n}d\left( X_{k},\bar{0}%
\right) \text{.}
\end{equation*}%
Finally, we have $\left( X_{k}\right) \in w^{\beta ,\infty }\left(
F,f\right) $ since $\beta \geq 1$ and $\left( X_{k}\right) \in w^{\beta
,\infty }\left( F\right) $ which completes the proof.

\bigskip

\noindent \textbf{Theorem 3.13.} Let $f$ be a modulus function such that $%
\lim \limits_{t\rightarrow \infty }\frac{f\left( t\right) }{t}>0$ and $\beta 
$ be a positive real number. Then $w^{\beta }\left( F,f\right) \subset
w^{\beta }\left( F\right) $.

\textbf{Proof.} Let $X=\left( X_{k}\right) $ be a sequence of fuzzy numbers
and $\left( X_{k}\right) \in w^{\beta }\left( F,f\right) .$ It is known that 
$\lim \limits_{t\rightarrow \infty }\frac{f\left( t\right) }{t}=\inf \left
\{ \frac{f\left( t\right) }{t}:t>0\right \} .$ We denote the value of $\lim
\limits_{t\rightarrow \infty }\frac{f\left( t\right) }{t}$ by $\ell $ for
shortness. Thus, we can write $f\left( t\right) \geq \ell t$ and $t\leq \ell
^{-1}f\left( t\right) $ for all $t\geq 0$ since $\ell >0$, and hence%
\begin{equation*}
\frac{1}{n^{\beta }}\dsum \limits_{k=1}^{n}d\left( X_{k},X_{0}\right) \leq
\ell ^{-1}\frac{1}{n^{\beta }}\dsum \limits_{k=1}^{n}f\left( d\left(
X_{k},X_{0}\right) \right) .
\end{equation*}%
Thus, we have $\left( X_{k}\right) \in w^{\beta }\left( F\right) .$

\bigskip

We have the following result from Theorem 3.12 and Theorem 3.13.

\bigskip

\noindent \textbf{Corollary 3.14.} Let $f$ be a modulus function. If $\lim
\limits_{t\rightarrow \infty }\frac{f\left( t\right) }{t}>0$ and $\beta \geq
1,$ then $w^{\beta }\left( F,f\right) =w^{\beta }\left( F\right) $.

\bigskip

\noindent \textbf{Theorem 3.15.} Let $f$ be a modulus function and $0<\beta
\leq \gamma $. Then $w^{\beta }\left( F,f\right) \subset w^{\gamma }\left(
F,f\right) $ and the inclusion is strict.

\textbf{Proof.} It is easy to show the inclusion relation. For strictness of
inclusion, let $f$ be a modulus function and consider the fuzzy sequence $%
X=\left( X_{k}\right) $ defined by%
\begin{equation*}
X_{k}\left( x\right) =\left \{ 
\begin{array}{cc}
\begin{array}{c}
x-1, \\ 
-x+3, \\ 
0%
\end{array}
& \left. 
\begin{array}{c}
1\leq x\leq 2 \\ 
2\leq x\leq 3 \\ 
\text{otherwise}%
\end{array}%
\right \} ,\text{ if }k=n^{3} \\ 
\bar{0}, & \text{if }k\neq n^{3}%
\end{array}%
\right. 
\end{equation*}%
We can write%
\begin{equation*}
\frac{1}{n^{\gamma }}\dsum \limits_{k=1}^{n}f\left( d\left( X_{k},\bar{0}%
\right) \right) \leq \frac{\sqrt[3]{n}}{n^{\gamma }}f\left( 3\right) =\frac{1%
}{n^{\gamma -\frac{1}{3}}}f\left( 3\right) 
\end{equation*}%
using the property $f\left( 0\right) =0$ for every $n\in 
\mathbb{N}
.$ So, we have $\left( X_{k}\right) \in w^{\gamma }\left( F,f\right) $ since
the right side tends to zero for $\gamma >\frac{1}{3}$ as $n\rightarrow
\infty .$ On the other hand, we obtain

\begin{equation*}
\frac{1}{n^{\beta }}\dsum \limits_{k=1}^{n}f\left( d\left( X_{k},\bar{0}%
\right) \right) \geq \frac{\sqrt[3]{n}-3}{n^{\beta }}f\left( 3\right)
\end{equation*}%
for every $n\in 
\mathbb{N}
.$ Hence we have $\left( X_{k}\right) \notin w^{\beta }\left( F,f\right) $
since the right side tends to infinity for $0<\gamma <\frac{1}{3}$ as $%
n\rightarrow \infty .$

\bigskip

\noindent \textbf{Theorem 3.16.} Let $X=\left( X_{k}\right) $ be a sequence
of fuzzy numbers and $0<\beta \leq \gamma \leq 1.$ Also, let $f$ be an
unbounded modulus such that there is a positive constant $c$ such that $%
f\left( xy\right) \geq cf\left( x\right) f\left( y\right) $ for all $x\geq
0, $ $y\geq 0$ and $\lim \limits_{t\rightarrow \infty }\frac{f\left(
t\right) }{t}>0$. If a sequence $X=\left( X_{k}\right) $ of fuzzy numbers is
strongly Cesaro summable of order $\beta $ with respect to modulus function $%
f$ to fuzzy number $X_{0},$ then it is $f-$statistically convergent of order 
$\gamma $ to fuzzy number $X_{0}.$

\textbf{Proof.} Take any fuzzy sequence $X=\left( X_{k}\right) $ and $%
\varepsilon >0.$ Then, we have%
\begin{eqnarray*}
\dsum \limits_{k=1}^{n}f\left( d\left( X_{k},X_{0}\right) \right) &\geq
&f\left( \dsum \limits_{k=1}^{n}d\left( X_{k},X_{0}\right) \right) \\
&\geq &f\left( \left \vert \left \{ k\leq n:d\left( X_{k},X_{0}\right) \geq
\varepsilon \right \} \right \vert \varepsilon \right) \\
&\geq &cf\left( \left \vert \left \{ k\leq n:d\left( X_{k},X_{0}\right) \geq
\varepsilon \right \} \right \vert \right) f\left( \varepsilon \right)
\end{eqnarray*}%
by the definition of modulus function, and so we can write 
\begin{eqnarray*}
\frac{1}{n^{\beta }}\dsum \limits_{k=1}^{n}f\left( d\left(
X_{k},X_{0}\right) \right) &\geq &\frac{cf\left( \left \vert \left \{ k\leq
n:d\left( X_{k},X_{0}\right) \geq \varepsilon \right \} \right \vert \right)
f\left( \varepsilon \right) }{n^{\beta }} \\
&\geq &\frac{cf\left( \left \vert \left \{ k\leq n:d\left(
X_{k},X_{0}\right) \geq \varepsilon \right \} \right \vert \right) f\left(
\varepsilon \right) }{n^{\gamma }} \\
&=&\frac{cf\left( \left \vert \left \{ k\leq n:d\left( X_{k},X_{0}\right)
\geq \varepsilon \right \} \right \vert \right) f\left( \varepsilon \right)
f\left( n^{\gamma }\right) }{n^{\gamma }f\left( n^{\gamma }\right) }
\end{eqnarray*}%
since $\beta \leq \gamma .$ Hence, it is easy to see that $X\in S^{\gamma
}\left( F,f\right) $ using the fact that $\lim \limits_{t\rightarrow \infty }%
\frac{f\left( t\right) }{t}>0$ and $X\in w^{\beta }\left( F,f\right) .$

\bigskip

We have the following result if we take $\beta =\gamma =1$ in Theorem 3.16.

\bigskip

\noindent \textbf{Corollary 3.17.} Let $f$ satisfies conditions in Theorem
3.16. If a sequence $X=\left( X_{k}\right) $ of fuzzy numbers is strongly
Cesaro summable with respect to modulus function $f$ to fuzzy number $X_{0},$
then it is $f-$statistically convergent to fuzzy number $X_{0}.$

\end{document}